\documentstyle[11pt,leqno]{article}
\newtheorem{theorem}{{\sc Theorem}}
\newcommand{\bt}{\begin{theorem}}
\newcommand{\et}{\end{theorem}}
\newcommand{\newsection}[1]{\setcounter{equation}{0} \setcounter{theorem}{0}
\section{#1}}

\newcommand{\NI}{\noindent}
\newcommand{\bea}{\begin{eqnarray}}
\newcommand{\eea}{\end{eqnarray}}

\def \spec#1 {\mathop{#1}}

\def \b #1 {\bf #1}

\newcommand{\ity}{\infty}
\newcommand{\raro}{\rightarrow}

\newcommand{\vsp}{\vskip 1em}

\newcommand{\be}{\begin{equation}}
\newcommand{\ee}{\end{equation}}
\newcommand{\ben}{\begin{eqnarray*}}
\newcommand{\een}{\end{eqnarray*}}

\begin{document}
\NI{\bf On some characterizations of probability distributions based on}\\
\NI{\bf maxima or minima of some families of dependent random variables}\\
\vsp
\NI{B.L.S. Prakasa Rao\footnote{CONTACT B.L.S. Prakasa Rao, blsprao@gmail.com, CR Rao advanced Institute of Mathematics, Statistics and Computer Science, Hyderabad 500046, India}}
\vsp
\NI{CR Rao Advanced Institute of Mathematics, Statistics and  Computer Science, Hyderabad, India}
\vsp
\NI{\bf Abstract}\\ 
Most of the characterizations of probability distributions are based on properties of functions of possibly independent random variables.  We investigate characterizations of  probability distributions through properties of minima or maxima of max-independent, min-independent and quasi-independent random variables generalizing the results from independent random variables of Kotlarski (1978), Prakasa Rao (1992) and Klebanov (1973).
\vsp
\NI{\bf Keywords}\\ 
Characterization; max-independence; min-independence; quasi-independence.
\vsp
\NI{\bf MSC 2020: Primary 62E10}\\
\vsp
\newsection{Introduction} Limit theorems as well as characterization results in probability and statistics are generally based on statistics which are functions of independent random variables. However  it was observed that the full force of independence of random variables is not necessary for justifying the conclusions. Durairajan (1979) introduced the concept of sub-independence for a set of random variables which is weaker than the notion of independent random variables. Hamedani (2013) gives a review of results in probability and statistics based on the concept of sub-independence of random variables. Hamedani and Prakasa Rao (2015) present characterization of probability distributions based on functions of random variables which are sub-independent or max-sub-independent or conditional-sub-independent. Kagan and Szekely (2016) introduced the notion of Q-independent random variables which is again an extension of the concept of independent random variables. Prakasa Rao (2018) proved an analogue of Skitovich-Darmois-Ramachandran-Ibragimov theorem for linear forms of Q-independent random variables and obtained a characterization of probability measures  based on linear functions of Q-independent random variables defined on a homogeneous Markov chain in Prakasa Rao (2022). There are several other results in literature based on Q-independent random variables which we are not mentioning here. We present characterizations for probability distributions through minima or maxima based on max-independent or min-independent or quasi-independent random variables, to be defined below, generalizing earlier results based on independent random variables.
\vsp
Let $X_0,X_1$ and $X_2$ be independent random variables. Define $Y_1= \max(X_0,X_1)$ and $Y_2= \max(X_0,X_2)$. It is of interest to know whether the joint distribution of $(Y_1,Y_2)$ determines the individual distributions of $X_0,X_1$ and $X_2$ uniquely. It is known that the random variable $Y_1$ alone can not determine the distributions of $X_0$ and $X_1$ uniquely unless $X_0$ and $X_1$ are identically distributed random variables (cf. Prakasa Rao (1992), Section 7.3).  Kotlarski (1978) and Klebanov (1973) obtained characterizations for probability distributions  through maxima or minima of independent random variables. We now discuss extensions of these results leading to characterizations of probability distributions through maxima and minima for some classes of dependent random variables.
\vsp
\newsection{Identifiability by maxima}
\vsp
The following result is due to Kotlarski (1978).
\vsp
\NI{\bf Theorem 2.1:} (Identifiability by maxima) Suppose $X_0,X_1,X_2$ are independent random variables. Define $Y_1= \max(X_0,X_1)$ and $Y_2= \max(X_0,X_2)$. Then the joint distribution of $(Y_1,Y_2)$ uniquely determines the distributions of the independent random variables $X_0,X_1$ and $X_2$ provided the supports of their distribution functions are the same.
\vsp
For a proof of Theorem 2.1, see Prakasa Rao (1992), Theorem 2.2.1, p.24.
\vsp
\NI {\bf Definition :} A finite collection of random variables $X_1,\dots,X_n$ is said to be {\it max-independent} if there exists a function $\eta(x_1,\dots,x_n)$ such that
$$F(x_1,\dots,x_n)=F_1(x_1)\dots F_n(x_n) \;\eta(x_1,\dots,x_n), x_i\in R, 1\leq i \leq n $$
where $ F(x_1,\dots,x_n)$ is the joint distribution of $(X_1,\dots,X_n), F_i(x)$ is the distribution function of $X_i$ for $1 \leq i \leq n$ and $\eta(x_1,\dots,x_n)$ is a  function taking values in the interval $(0,1]$ such that $\eta(x_1,\dots,x_n) \raro 1$ if $x_i \raro \ity$ for some $i, 1\leq i \leq n.$ It is said to be {\it min-independent} if there exists a function $\eta(x_1,\dots,x_n)$ such that
$$P(X_1>x_1,\dots,X_n>x_n)= P(X_1>x_1)\dots P(X_n>x_n) \;\eta(x_1,\dots,x_n), x_i\in R, 1\leq i \leq n $$
where $\eta(x_1,\dots,x_n)$ is a  function taking values in the interval $(0,1]$ such that $\eta(x_1,\dots,x_n) \raro 0$ if $x_i \raro -\ity$for some $i, 1\leq i \leq n.$ The sequence of random variables $X_1,\dots,X_n$ is said to be {\it quasi-independent} if there exists a  measure $\eta(.) \neq 0$ on $R^n$ such that $\eta(B_1 X \dots X B_n)=1$ whenever at least one of the sets $B_i$ is $R$ and 
$$P(X_1\in B_1,\dots, X_n\in B_n)= \Pi_{j=1}^n P(X_j\in B_j)\;\eta(B_1 X\dots X B_n)$$
for all Borel subsets $B_i, 1\leq i \leq n$ in $R.$ Hereafter we call the measure $\eta(.)$ as the {\it generator} of the random vector $(X_1,\dots,X_n).$
\vsp
Let $F_i(.)$ be the distribution function of the random variable $X_i$ for $1\leq i \leq n$ in the following discussion. It is easy to see that every set of independent random variables $X_1,\dots,X_n$ is max-independent but a set of max-independent random variables need not be independent. For instance, the random variables $Y_1$ and $Y_2$ defined above are not independent but they are max-independent since
\ben
F_{Y_1,Y_2}(y_1,y_2)&=&P(Y_1 \leq y_1, Y_2\leq y_2)\\
&=& P(\max(X_0,X_1) \leq y_1, \max(X_0,X_2) \leq y_2)\\
&=& P(X_0 \leq y_1, X_1 \leq y_1, X_0 \leq y_2, X_2 \leq y_2)\\
&=& P(X_0 \leq \min(y_1,y_2), X_1 \leq y_1, X_2\leq y_2)\\
&=& F_0(\min(y_1,y_2))F_1(y_1)F_2(y_2)\\
&=& \frac{F_0(\min(y_1,y_2))}{F_0(y_1)F_0(y_2)}F_0(y_1)F_0(y_2)F_1(y_1)F_2(y_2)\\
&=& F_{Y_1}(y_1)F_{Y_2}(y_2) \eta(y_1,y_2)\\
\een
where
$$\eta(y_1,y_2)= \frac{F_0(\min(y_1,y_2))}{F_0(y_1)F_0(y_2)}.$$ 
Note that 
$$F_{Y_1}(y_1)= F_0(y_1)F_1(y_1)$$
and
$$ F_{Y_2}(y_2)= F_0(y_2)F_2(y_2).$$
Here $F_{Y}(.)$ denotes the distribution function of the random variable $Y$. Observe that 
$\eta(y_1,y_2) \raro 1$ as $y_1\raro \ity$ or $y_2 \raro \ity$ from the properties of a distribution function. Hence the random variables $Y_1$ and $Y_2$ are max-independent. We will now prove an extension of Theorem 2.1 leading to identifiability by maxima.
\vsp
\NI{\bf Theorem 2.2:} (Identifiability by maxima) Let $X_0,X_1$ and $X_2$ be max-independent random variables. Define $Y_1=\max(X_0,X_1)$ and $Y_2=\max(X_0,X_2).$ Then the  joint distribution of the bivariate random vector $(Y_1,Y_2)$ uniquely determines the distributions of $X_0,X_1$ and $X_2$ provided the supports of the distribution functions of $X_0,X_1$ and $X_2$ are the same.
\vsp
\NI{\bf Proof:} Let $F_0(.),F_1(.) $ and $F_2(.)$  be the distribution functions of $X_0,X_1$ and $X_2$ respectively and $G(y_1,y_2)$ be the joint distribution of the bivariate random vector $(Y_1,Y_2).$ Observe that
\ben
G(y_1,y_2) &= & P(Y_1\leq y_1,Y_2\leq y_2)\\
&=& P( X_0\leq y_1, X_1\leq y_1; X_0\leq y_2, X_2\leq y_2)\\
&=& P(X_0\leq \min(y_1,y_2),X_1\leq y_1,X_2\leq y_2)\\
&= & F_0(\min(y_1,y_2))F_1(y_1)F_2(y_2) \eta_1(\min(y_1,y_2),y_1,y_2)\\
\een
by the max-independence of the random variables $X_0,X_1,X_2$ for some positive function $\eta_1(z,y_1,y_2)$  which tends to one as $z \raro \ity$ or $y_1 \raro \ity$ or  $y_2\raro \ity.$ Let $F_0^*(.),F_1^*(.) $ and $F_2^*(.)$  be the alternate distribution functions of $X_0^*,
X_1^*$ and $X_2^*$ respectively which are max-independent. Define $Y_1^*=\max(X_0^*,X_1^*)$ ans $Y_2^*=\max(X_0^*, X_2^*).$  Suppose that the joint distribution of $(Y_1^*,Y_2^*)$ is the same as that of $(Y_1,Y_2).$  Following calculations similar to those given above, it follows that
\be
G(y_1,y_2)=F_0^*(\min(y_1,y_2))F_1^*(y_1)F_2^*(y_2) \eta_2(\min(y_1,y_2),y_1,y_2)
\ee
for for some non-negative function $\eta_2(z,y_1,y_2)$  which tends to one as $z \raro \ity$ or $y_1 \raro \ity$ or $y_2\raro \ity.$ Hence
\begin{eqnarray*}
F_0(\min(y_1,y_2))F_1(y_1)F_2(y_2) \eta_1(\min(y_1,y_2),y_1,y_2)\\
=F_0^*(\min(y_1,y_2))F_1^*(y_1)F_2^*(y_2) \eta_2(\min(y_1,y_2),y_1,y_2), y_1,y_2\in R.\\
\end{eqnarray*} 
Fix $y_1\leq y_2$. It follows that
\be
F_0(y_1)F_1(y_1)F_2(y_2) \eta_1(y_1,y_1,y_2)=F_0^*(y_1)F_1^*(y_1)F_2^*(y_2) \eta_2(y_1,y_1,y_2), y_2\in R.
\ee
Let $y_2 \raro \ity.$ Then it follows that
\be
F_0(y_1)F_1(y_1)= F_0^*(y_1)F_1^*(y_1), y_1\in R 
\ee
by the properties of distribution functions. Equations (2.2) and (2.3) show that
$$ F_2(y_2)\eta_1(y_1,y_1,y_2)=F_2^*(y_2) \eta_2(y_1,y_1,y_2), y_2\in R$$
provided $F_0(y_1)F_1(y_1) >0.$ Let $y_1\raro \ity.$ It follows that
\be
F_2(y_2)=F_2^*(y_2), y_2\in R .
\ee
Note that the support of the function $F_0F_1$ is the same as the support of the function $F_0^*F_1^*$ from (2.3). Similar analysis proves that
\be
F_1(y_1)=F_1^*(y_1), y_1\in R
\ee
provided $F_0(y_2)F_2(y_2)>0.$Note that the support of $F_0F_2$ is the same as that of the support of $F_0^*F_2^*.$ Since the supports of $F_0,F_1$ and $F_2$ are all the same, it follows from (2.2)-(2.5) that
$$F_i(y)=F_i^*(y), i=0,1,2$$
over the common support of the distribution functions of $X_0,X_1$ and $X_2.$ Hence the distribution of the bivariate random vector $(Y_1,
Y_2)$ uniquely determines the distributions of $X_0,X_1$ and $X_2.$ 
\vsp
\newsection{Identifiability by minima}
\vsp
The following result is known for the minima of random variables.
\vsp
\NI{\bf Theorem 3.1:} (Identifiability by minima) Suppose $X_0,X_1,X_2$ are independent random variables. Define $Y_1= \min(X_0,X_1)$ and $Y_2= \min(X_0,X_2)$. Suppose the distribution functions $F_0,F_1$ and $F_2$ of $X_0,X_1$ and $X_2$ respectively satisfy the conditions
$$F_i(a)=1, F_i(x)<1, x<a, i=0,1,2$$
for some $a \leq \infty.$ Then the joint distribution of $(Y_1,Y_2)$ uniquely determines the distributions of the independent random variables $X_0,X_1$ and $X_2$.
\vsp
For a proof of Theorem 3.1, see Kotlarski (1978) and Prakasa Rao (1992), Theorem 2.3.1, p.24.
\vsp
We now obtain an extension of Theorem 3.1 for min-independent random variables.
\vsp
\NI{\bf Theorem 3.2:}  (Identifiabilty by minima) Let $X_0,X_1$ and $X_2$ be min-independent random variables. Define $Y_1=\min(X_0,X_1)$ and $Y_2=\min(X_0,X_2).$ Suppose the distribution functions $F_0,F_1$ and $F_2$ of $X_0,X_1$ and $X_2$ respectively satisfy the conditions
$$F_i(a)=1, F_i(x)<1, x<a, i=0,1,2$$
for some $a \leq \infty.$ Then the  joint distribution of the bivariate random vector $(Y_1,Y_2)$ uniquely determines the distributions of $X_0,X_1$ and $X_2$.
\vsp
\NI{\bf Proof:} Let $\bar F_i=1-F_i, i=0,1,2$ and $S(y_1,y_2)=P(Y_1>y_1,Y_2>y_2), y_1,y_2\in R.$ Then,  
\ben
S(y_1,y_2) &= & P(Y_1 > y_1,Y_2 > y_2)\\
&=& P( X_0 > y_1, X_1 >y_1; X_0 > y_2, X_2 > y_2)\\
&=& P(X_0 > \max(y_1,y_2),X_1 > y_1,X_2 >y_2)\\
&= & \bar F_0(\max(y_1,y_2))\bar F_1(y_1) \bar F_2(y_2) \eta_1(\max(y_1,y_2),y_1,y_2).\\
\een
Let $F_0^*(.),F_1^*(.) $ and $F_0^*(2)$  be the alternate distribution functions of $X_0,X_1$ and $X_2$ respectively which are min-independent  such that the joint distribution of $(Y_1^*,Y_2^*)$ is the same as that of $(Y_1,Y_2)$ where 
$Y_1^*=\min(X_0^*,X_1^*)$ and $Y_2^*=\min(X_0^*,X_2^*).$ Following similar calculations as given above, it follows that
\be
S(y_1,y_2)=\bar F_0^*(\max(y_1,y_2))F_1^*(y_1)F_2^*(y_2) \eta_2(\max(y_1,y_2),y_1,y_2)
\ee
for for some non-negative function $\eta_2( \max(y_1,y_2),y_1,y_2)$  which tends to zero as $y_1 \raro -\ity$ or $y_2\raro -\ity.$
Following arguments similar to those given in the proof of Theorem 2.2, it follows that $\bar F_i=\bar F_i^*, i=0,1,2$ and hence the  joint distribution of the bivariate random vector $(Y_1,Y_2)$ uniquely determines the distributions of $X_0,X_1$ and $X_2$ provided there exists some $a \leq \infty$ such that 
$$F_i(a)=F_i^*(a)=1; F_i(x)<1, x<a;F_i^*(x)<1, x<a, i=0,1,2.$$
\vsp
\newsection{Identifiability by maximum and minimum}
\vsp
Let $X_0,X_1$ and $X_2$ be independent random variables. Define $Y_1= \min(X_0,X_1)$ and $Y_2=\max(X_0,X_2).$ The following theorem, due to Kotlarski,  holds leading to identifiability of the distributions of $X_0,X_1$ and $X_2$ given the distribution of the bivariate random vector $(Y_1,Y_2).$ 
\vsp
\NI{\bf Theorem 4.1:} (Identifiability by maximum and minimum) Let $X_0,X_1$ and $X_2$ be independent random variables and $F_i$ be the distribution function of $X_i, i=0,1,2.$ Suppose that there exists $a,b,x_0,q$ satisfying $-\infty \leq a,x_0< b\leq \infty, 0 < q <1,$ such that
$$F_1(x)<1, x<b; F_1(b-0)=1 \;\;\mbox{if}\;\;  b \in R,$$
$$F_2(y)>0, y>a; F_2(a+0)=1 \;\;\mbox{if}\;\;  a \in R,$$
$$F_0(a+0)=0,F_0(b-0)=1,F_0(x_0)=q$$
and $F_0$ is strictly increasing in the interval $(a,b).$ Then the joint distribution of $(Y_1,Y_2)$ where  $Y_1= \min(X_0,X_1)$ and $Y_2=\max(X_0,X_2)$ uniquely determines the distributions of $F_0,F_1$ and $F_2.$
\vsp
For the proof of Theorem 4.1, see Kotlarski (1978) or Prakasa Rao (1992),Theorem 2.4.1. 
\vsp
We will now extend the result in Theorem 4.1 to quasi-independent random variables $X_i, i=0,1,2.$
\vsp
\NI{\bf Theorem 4.2:} (Identifiability by maximum and minimum) Let $X_0,X_1$ and $X_2$ be quasi-independent random variables with a given generator $\eta(.)$ and $F_i$ be the distribution function of $X_i, i=0,1,2.$ Suppose that there exists $a,b,x_0,q$ satisfying $-\infty \leq a,x_0< b\leq \infty, 0 < q <1,$ such that
$$F_1(x)<1, x<b; F_1(b-0)=1 \;\;\mbox{if}\;\;  b \in R,$$
$$F_2(y)>0, y>a; F_2(a+0)=1 \;\;\mbox{if}\;\; a \in R,$$
$$F_0(a+0)=0,F_0(b-0)=1,F_0(x_0)=q$$
and $F_0$ is strictly increasing in the interval $(a,b).$ Then the joint distribution of $(Y_1,Y_2),$ where $Y_1= \min(X_0,X_1)$ and $Y_2=\max(X_0,X_2),$ uniquely determines the distributions of $F_0,F_1$ and $F_2$ given the generator $\eta(.).$
\vsp
\NI{\bf Proof:} Let $-\infty <y_1 < y_2<\infty.$ Then
\begin{eqnarray*}
P(Y_1>y_1,Y_2\leq y_2)&=& P(X_0>y_1,X_1>y_1, X_0\leq y_2,X_2\leq y_2)\\
&=& P( y_1 < X_0 \leq y_2, X_1>y_1, X_2\leq y_2)\\
&=& P(X_0 \in B_0, X_1\in B_1, X_2\in B_2)\\
&=& P(X_0 \in B_0) P(X_1\in B_1) P(X_2\in B_2) \eta(B_0 X B_1 X B_2)\\
\end{eqnarray*}
where $B_0= (y_1,y_2], B_1= (y_1,\infty), B_2=(-\infty, y_2]$ for some measure $\eta(.)$ on $R^3$ with properties as mentioned earlier by the property of quasi-independence of the random variables $X_0,X_1$ and $X_2.$ It is easy to see that
\ben
P(Y_1>y_1,Y_2\leq y_2)= [F_0(y_2)-F_0(y_1)]\bar F_1(y_1)F_2(y_2) \eta(B_0 X B_1 X B_2), y_1,y_2\in R
\een
where $\bar F(x)= 1-F(x).$ Suppose $\{F_0^*,F_1^*,F_2^*\}$ ia another set of distributions functions for $\{X_0,X_1,X_2\}$ with the generator $\eta(.)$ satisfying the conditions stated in the theorem such that the distribution functions of $(Y_1,Y_2)$ under $\{F_0^*,F_1^*,F_2^*\}$ as well as $\{F_0,F_1,F_2\}$ are the same. Then, for 
$-\infty <y_1 \leq y_2<\infty,$
\be
[F_0(y_2)-F_0(y_1)]\bar F_1(y_1)F_2(y_2) \eta(B_0 X B_1 X B_2)=[F_0^*(y_2)-F_0^*(y_1)]\bar F_1^*(y_1)F_2^*(y_2) \eta(B_0 X B_1 X B_2)
\ee
with properties of the measure $\eta(.)$ as mentioned earlier in the definition of quasi-independence.  Let $y_2 \raro \ity$ in (4.1). Since $B_2 \uparrow R$ as $y_2\raro \infty$, it follows that
\be
\bar F_0(y_1) \bar F_1(y_1)= \bar F_0^*(y_1) \bar F_1^*(y_1), y_1\in R
\ee
from the properties of the distribution functions. Let $y_1\raro -\infty$ in (4.1). Then $B_1\uparrow R$ as $y_1\raro -\infty$ and it follows again from (4.1) that
\be
F_0(y_2) F_2(y_2)= F_0^*(y_2) F_2^*(y_2), y_2\in R.
\ee
Combining the equations (4.1)-(4.3), it follows that
\begin{eqnarray*}
\lefteqn{[F_0(y_2)-F_0(y_1)]\bar F_1(y_1)F_2(y_2) \eta(B_0 X B_1 X B_2) \bar F_0^*(y_1) \bar F_1^*(y_1)\bar F_0^*(y_2) \bar F_2^*(y_2)}\\
&=&[F_0^*(y_2)-F_0^*(y_1)]\bar F_1^*(y_1)F_2^*(y_2) \eta(B_0 X B_1 X B_2) \bar F_0(y_1) \bar F_1(y_1)\bar F_0(y_2) \bar F_2(y_2) 
\end{eqnarray*}
for $-\infty<y_1 < y_2<\infty.$ Applying the conditions stated in the theorem connected with the supports of the distribution functions $F_i, i=0,1,2$ and of $F_i^*, i=0,1,2,$ it follows that
\be
\frac{F_0^*(y_2)-F_0^*(y_1)}{F_0(y_2)-F_0(y_1)}= \frac{\bar F_0^*(y_1)}{\bar F_0(y_1)}\frac{F_0^*(y_2)}{F_0(y_2)}, -\infty\leq a<y_1<y_2<b\leq \infty. 
\ee
Let $y_1=y$ and $y_2=x_0$ in the equation given above. Since $F_0^*(x_0)=F_0(x_0)=q,$ it follows that 
\be
\frac{F_0^*(x_0)-F_0^*(y)}{F_0(x_0)-F_0(y)}= \frac{\bar F_0^*(y)}{\bar F_0(y)}, -\infty \leq a<y < x_0. 
\ee
This relation in turn implies that
\be
F_0^*(y)=F_0(y), -\infty <y \leq x_0 
\ee 
after a simple algebra. Similarly we can prove that
\be
F_0^*(y)=F_0(y), x_0\leq y<\infty.
\ee 
Equations (4.6) and (4.7) show that $F_0^*(y)=F_0(y), -\infty<y<\infty$. It then follows that $F_1^*(y)=F_1(y), -\infty<y<\infty$ and $F_2^*(y)=F_2(y), -\infty<y<\infty$ from the equations (4.2) and (4.3) completing the proof of the theorem.
\vsp
\newsection{Identifiability by maxima of several random variables}
\vsp
The following result, due to Klebanov (1973), deals with identifiability by maxima of several random variables.
\vsp
\NI{\bf Theorem 5.1:} Suppose that $X_1,\dots,X_n$ are independent positive random variables with distribution functions $F_1,\dots,F_n$ respectively. Further suppose that $F_i(x)>0$ for $x>0, 1\leq i \leq n.$ Define
$$Y_1=\max(a_1X_1,\dots,a_n X_n)\;\;\mbox{and}\;\; Y_2=\max(b_1X_1,\dots,b_n X_n)$$
where $a_i>0,b_i>0, 1\leq i \leq n$ and $a_i:b_i\neq a_j:b_j, 1\leq i\neq j\leq n.$ Then the joint distribution of $(Y_1,Y_2)$  uniquely determines the distribution functions of $X_i, 1\leq i \leq n.$
\vsp
For proof of this result, see Prakasa Rao (1992), Theorem 2.8.1. 
\vsp
We will now extend this result to random vector $(X_1,\dots,X_n)$ where $X_1,\dots,X_n$ are positive max-independent random variables.
\vsp
\NI{\bf Theorem 5.2:} Suppose that $X_1,\dots,X_n$ are max-independent positive random variables with distribution functions $F_1,\dots,F_n$ respectively with a given generator $\eta(x_1,\dots,x_n)).$ Further suppose that $F_i(x)>0$ for $x>0, 1\leq i \leq n.$ Define
$$Y_1=\max(a_1X_1,\dots,a_nX_n)\;\; {\mbox and}\;\; Y_2=\max(b_1X_1,\dots,b_nX_n)$$
where $a_i>0,b_i>0, 1\leq i \leq n$ and $a_i:b_i\neq a_j:b_j, 1\leq i\neq j\leq n.$ Then the joint distribution of $(Y_1,Y_2)$  uniquely determines the distribution functions of $X_i, 1\leq i \leq n.$
\vsp
\NI{\bf Proof:} Let $G(t,s)$ be the joint distribution function of the bivariate random vector $(Y_1,Y_2)$ as defined above. Note that
\begin{eqnarray*}
P(Y_1\leq t, Y_2\leq s)&=& P(\max(a_1X_1,\dots,a_nX_n)\leq t; \max(b_1X_1,\dots,b_nX_n)\leq s)\\
&=& P(X_i\leq \min(t/a_i;s/b_i), 1\leq i \leq n)\\
&=& \Pi_{i=1}^n P(X_i\leq \min(t/a_i;s/b_i)) \eta(\min{t/a_i;s/b_i}, 1\leq i \leq n)\\
&=& \Pi_{i=1}^nF_i(\min(t/a_i;s/b_i))\eta(\min(t/a_i;s/b_i), 1\leq i \leq n)\\
\end{eqnarray*}
for $0\leq t,s<\infty.$ Suppose $F_i*, 1\leq i \leq n$ alternate distributions of $X_i, 1\leq i \leq n$ satisfying the conditions in the theorem. Then, it follows that
\begin{eqnarray*}
\Pi_{i=1}^nF_i(\min(t/a_i;s/b_i))\eta(\min(t/a_i;s/b_i), 1\leq i \leq n)\\
=\Pi_{i=1}^nF_i^*(\min(t/a_i;s/b_i))\eta(\min(t/a_i;s/b_i), 1\leq i \leq n), 0\leq t,s<\infty
\end{eqnarray*}
or equivalently
\be
\Pi_{i=1}^nF_i(\min(t/a_i;s/b_i))=\Pi_{i=1}^nF_i^*(\min(t/a_i;s/b_i)), 0\leq t,s<\infty.
\ee
Let $v_j(t)=\log F_i(t/b_i)-\log F_i^*(t/b_i).$ Then the equation (5.2) can be written in the form
\be
\sum_{i=1}^nv_i(\min(c_it,s))=0, 0\leq t,s<\infty
\ee
where $c_i=\frac{b_i}{a_i}, 1\leq i \leq n$ are pairwise distinct nonzero positive numbers. Without loss of generality, we assume that $0<
c_1<\dots<c_n<\infty.$ Let $t>0$ and $s=\tau t$ where $c_{n-1}<\tau<c_n.$ Then the equation (5.2) can be written in the form
\ben
\sum_{i=1}^{n-1}v_i(c_it)+ v_n(\tau t)=0, 0<t<\infty.
\een
This equation shows that $v_n(.)$ is constant on the interval $(c_{n-1}t,c_nt)$ for any $t>0.$ Since $t>0$ is arbitrary, it follows that $v_n(.)$ is constant on the interval $(0,\infty).$ Since $v_i(t) \raro 0$ as $t\raro \ity$ from the properties of distribution functions, it follows that $v_n(t)=0$ for $t>0.$ Repeating this process, it follows that $v_i(t)=0, 1\leq i \leq n-1, t>0.$ This, in turn implies that
\be
F_i(t/b_i)=F_i^*(t/b_i), 1\leq i\leq n, 0<t<\infty
\ee
from the definition of the function $v_i(t).$ Since $t>0$ is arbitrary, it follows that
\be
F_i(t)=F_i^*(t), 1\leq i \leq n, 0<t<\infty
\ee
proving the theorem.
\vsp
We now extend Theorem 5.1 to max-independent random variables which are not necessarily positive under some conditions.
\vsp
\NI{\bf Theorem 5.3:} Suppose that $X_1,\dots,X_n$ are max-independent random variables with distribution functions $F_i(x)>0, x \in R$ and $P(X_i=0)=0$ for  $1\leq i \leq n$ and with a given generator $\eta(x_1,\dots,x_n)).$ Define
$$Y_1=\max(a_1X_1,\dots,a_nX_n)\;\; {\mbox and}\;\; Y_2=\max(b_1X_1,\dots,b_nX_n)$$
where $a_i>0,b_i>0, 1\leq i \leq n$ and $a_i:b_i\neq a_j:b_j, 1\leq i\neq j\leq n.$ Then the joint distribution of $(Y_1,Y_2)$  uniquely determines the distribution functions of $X_i, 1\leq i \leq n.$
\vsp
\NI{\bf Proof:} As described in the proof of Theorem 5.2, it follows that 
\be
\sum_{i=1}^n v_i(\min(c_it,s))=0, 0\leq t,s<\infty
\ee
where $c_i=\frac{b_i}{a_i}, 1\leq i \leq n$ are pairwise distinct nonzero positive numbers. Without loss of generality, we assume that $0<
c_1<\dots<c_n<\infty.$ Following the arguments given in the proof of Theorem 5.2, it follows that $v_i(t)=0,t>0.$ Suppose $t<0.$ Let $s=\tau t, \tau \in (c_1,c_2).$ Then the equation takes the form
\be
v_1(\tau t)+\sum_{i=2}^n v_i(c_it)=0.
\ee
Hence $v_1(.)$ is constant on the interval $(c_2t,c_1t).$ since $t<0$ is arbitrary, it follows that $v_(t)=0$ on $(-\infty,0)$. Note that $v_1(t)$ is continuous at $t=0$ Hence $v_1(0)=0.$ Therefore $v_1(t)=0$ for all $t\in R.$ Applying induction arguments, it follows that $v_i(t)=0,$
for all $t\in R, 1\leq j \leq n.$ Hence $F_i=F_i^*, 1\leq i \leq n$ completing the proof of the theorem.
\vsp
\NI{\bf Remarks:} Following the arguments similar to those given in Theorem 5.3, one can prove an analogous result for min-independent random variables leading to the following theorem. we omit the details.
\vsp
\NI{\bf Theorem 5.4:} Suppose that $X_1,\dots,X_n$ are min-independent random variables with distribution functions $F_i(x)>0, x \in R$ and $P(X_i=0)=0$ for  $1\leq i \leq n$ and with a given generator $\eta(x_1,\dots,x_n)).$ Define
$$Y_1=\min(a_1X_1,\dots,a_nX_n) ; Y_2=\min(b_1X_1,\dots,b_nX_n)$$
where $a_i>0,b_i>0, 1\leq i \leq n$ and $a_i:b_i\neq a_j:b_j, 1\leq i\neq j\leq n.$ Then the joint distribution of $(Y_1,Y_2)$  uniquely determines the distribution functions of $X_i, 1\leq i \leq n.$
\vsp
\NI{\bf Acknowledgment:} This work was supported under the scheme ``INSA Honorary Scientist" at the CR Rao Advanced Institute of Mathematics, Statistics and Computer Science, Hyderabad 500046, India.
\vsp
\NI{\bf References}
\begin{description}
\item Durairajan, T.M. 1979. A class room mnote on sub-independence, {\it Gujarat Stat. Rev.}, 6:17-18.

\item Hamedani, G.G. 2013. Sub-independence: An expository perspective, {\it Commun. Statist. Theory-Methods}, 41:3615-3638.

\item Hamedani, G.G. and Prakasa Rao, B.L.S. 2015. Characterizations of probability distributions thriugh sub-independence, max-sub-independence and conditional sub-independence, {\it J. Indian Stat. Assoc.}, 53:79-87.

\item Kagan, A.M. and G. J. Szekely. 2016. An analytic generalization of independent and identical distributiveness, {\it Statist. Probab. Lett.}, 110:244-248.

\item Klebanov, L. 1973. Reconstructing the distributions of the components of a random vector from distributions of certain statistics, {\it Mathematical Notes}, 13:71-72.

\item Kotlarski, I. 1978. On some characterization in probability by using minima and maxima of random variables, {\it Aequationes Mathematicae}, 17:77-82.

\item Prakasa Rao, B.L.S. 1992. {\it Identifiability in Stochastic Models: Characterization of Probability Distributions}. New York: Academic Press.

\item Prakasa Rao, B.L.S. 2018. On the Skitovich-Darmois-Ramachandran-Ibragimov theorem for linear forms of q-independent random variables, {\it Studia  Scientarium Mathematicarum Hungarica}, 55:353-363.

\item Prakasa Rao, B.L.S. 2022. Characterization of probability measures by linear functions of Q-independent random variables defined on a homogeneous Markov chain, {\it Commun. Statist. Theory-Methods}, 51:6529-6534.
\end{description}

\end{document}